\documentclass[11.05pt,a4paper]{amsart}
\usepackage{amssymb,amsmath}
\usepackage{color}
\usepackage{latexsym}
\usepackage{amsthm,amsfonts,amssymb,mathrsfs}
\usepackage{rotating}
\usepackage[leqno]{amsmath}
\usepackage{xspace}
\usepackage[all]{xy}
\usepackage{longtable}

\headsep=1cm \numberwithin{equation}{section}
\newtheorem{theorem}{Theorem}[section]
\newtheorem{lemma}[theorem]{Lemma}
\newtheorem{proposition}[theorem]{Proposition}
\newtheorem{corollary}[theorem]{Corollary}

\theoremstyle{definition}
\newtheorem{definition}[theorem]{Definition}
\theoremstyle{remark}
\newtheorem{remark}[theorem]{Remark}

\newtheorem{acknowledgement}{Acknowledgement}

\newcommand{\fdepth}{\operatorname{fdepth}}

\newcommand{\Spec}{\operatorname{Spec}}

\newcommand{\amp}{\operatorname{amp}}

\newcommand{\cd}{\operatorname{cd}}

\newcommand{\Ht}{\operatorname{ht}}
\newcommand{\id}{\operatorname{id}}
\newcommand{\fd}{\operatorname{fd}}
\newcommand{\Gid}{\operatorname{Gid}}

\newcommand{\Gpd}{\operatorname{Gpd}}
\newcommand{\V}{\operatorname{V}}

\newcommand{\Cone}{\operatorname{Cone}}
\newcommand{\Tel}{\operatorname{Tel}}
\newcommand{\Ext}{\operatorname{Ext}}
\newcommand{\Supp}{\operatorname{Supp}}
\newcommand{\Tor}{\operatorname{Tor}}
\newcommand{\Hom}{\operatorname{Hom}}

\newcommand{\depth}{\operatorname{depth}}
\newcommand{\width}{\operatorname{width}}

\newcommand{\vpl}{\operatornamewithlimits{\varprojlim}}

\newcommand{\lo}{\longrightarrow}
\newcommand{\fm}{\frak{m}}
\newcommand{\fp}{\frak{p}}

\newcommand{\fa}{\frak{a}}
\newcommand{\fb}{\frak{b}}

\newenvironment{prf}[1][Proof]{\begin{proof}[\bf #1]}{\end{proof}}
\begin{document}

\author[Asgharzadeh and Divaani-Aazar]{Mohsen Asgharzadeh and
Kamran Divaani-Aazar$^\ast$ }

\title[A unified approach to formal local ...]
{A unified approach to formal local cohomology and local Tate
cohomology}

\address{M. Asgharzadeh, School of Mathematics, Institute for
Research in Fundamental Sciences (IPM), P.O. Box 19395-5746, Tehran,
Iran.} \email{asgharzadeh@ipm.ir}
\address{K. Divaani-Aazar, Department of Mathematics, Az-Zahra
University, Vanak, Post Code 19834, Tehran, Iran-and-School of
Mathematics, Institute for Research in Fundamental Sciences (IPM),
P.O. Box 19395-5746, Tehran, Iran.} \email{kdivaani@ipm.ir}


\keywords{Cohen-Macaulay complexes; formal depth; formal local
cohomology; Frobenius depth; local cohomology;
local homology; local Tate cohomology.\\
The second author was supported by a grant from IPM (No. 88130122).\\
$^\ast$ Corresponding author}

\begin{abstract} Let $R$ be a commutative Noetherian ring.
We introduce a theory of formal local cohomology for complexes of
$R$-modules. As an application, we establish some relations between
formal local cohomology, local homology, local cohomology and
local Tate cohomology through some natural isomorphisms. We
investigate vanishing of formal local cohomology modules. Also,
we give a characterization of Cohen-Macaulay complexes.
\end{abstract}

\maketitle

\section{Introduction}

Throughout,  $R$ is a commutative Noetherian ring with a nonzero
identity. Let $\fa$ and $\fb$ be two ideals of $R$ and $M$ a
finitely generated $R$-module. When $R$ is local with maximal ideal
$\fm$, for each $i\geq 0$, Schenzel \cite{Sc2} has called
$\mathfrak{F}^i_{\fa}(M)={\vpl}_nH^i_{\fm}(M/\fa^nM)$ $i$-th formal
local cohomology module of $M$ with respect to $\fa$. Our aim in
this note is to establish a theory of formal local cohomology in
$\mathcal{D}(R)$, the derived category of $R$-modules. Let ${\bf
R}\Gamma_{\fa}(-)$ and ${\bf L}\Lambda^{\fa}(-)$ denote derived
local cohomology and derived local homology functors with respect to
$\fa$. (We recall their definitions in the beginning of the next section.)
The compositions of these derived functors were studied extensively
in \cite{LLT} and
\cite{L}. For a complex $X\in
\mathcal{D}(R)$, we call $\mathfrak{F}_{\fa,\fb}(X):={\bf
L}\Lambda^{\fa}({\bf R}\Gamma_{\fb}(X))$ formal local complex of $X$
with respect to $(\fa,\fb)$. Also, for each integer $i$, we call
$\mathfrak{F}_{\fa,\fb}^i(X):= H_{-i}(\mathfrak{F}_{\fa,\fb}(X))$
$i$-th formal local cohomology module of $X$ with respect to
$(\fa,\fb)$. It is worth to point out that in the case $R$ is local
with maximal ideal $\fm$, we show that
$\mathfrak{F}_{\fa,\fm}^i(M)=\mathfrak{F}_{\fa}^i(M)$ for all $i$,
see Corollary 2.4 below.

In Section 2, we establish several general properties of formal
local cohomology. In the case that $R$ possesses a normalized
dualizing complex, we present a duality result, see Lemma 2.2 below.
This duality lemma facilitate working with formal local cohomology
modules. We continue this section by applying the theory of formal
local cohomology for examining local Tate cohomology modules. Let
$X\in \mathcal{D}(R)$. The local Tate cohomology modules
$\widehat{H}^i_{\fa}(X)$ were introduced by Greenlees \cite{G} and
their study was continued by  Lipman,  L\'{o}pez and Tarrioin
\cite{LLT}. We deduce a long exact sequence
$$\cdots \lo H_{\fa}^i(X)\lo \mathfrak{F}_{\fa,\fa}^i(X)\lo
\widehat{H}^i_{\fa}(X)\lo H_{\fa}^{i+1}(X)\lo
\mathfrak{F}_{\fa,\fa}^{i+1}(X)\lo \cdots, $$ which, in turn, yields
several corollaries. In particular, from this sequence, we
immediately deduce that  $\widehat{H}_{\fa}^i(M)\cong
H_{\fa}^{i+1}(M)$ for all $R$-modules $M$ and all $i>0$.

Section 3 is the core of this paper. In this section,  we present
several results on vanishing of formal local cohomology modules. Let
$X$ be a complex. We introduce notion $\fdepth(\fa,\fb,X)$ as the
infimum of the integers $i$ such that
$\mathfrak{F}_{\fa,\fb}^i(X)\neq 0$. In the prime characteristic
case, this notion is closely related to the notion of Frobenius
depth which was defined by Hartshorne and Speiser in \cite{HS}. When
$X$ is homologically bounded and all of its homology modules are
finitely generated, we establish the inequality
$$\depth(\fb,M)-\cd_{\fa}(R)\leq \fdepth(\fa,\fb,X)\leq \dim_RX/\fa
X.$$ Let $(R,\fm)$ be a local ring and $0\not\simeq X$ a homologically
bounded complex with finitely generated homology modules. We show
that the supremum of the integers $i$ such that
$\mathfrak{F}_{\fa,\fm}^i(X)\neq 0$ is equal to $\dim_RX/\fa X$. Set
$H(X)^\sharp:=\oplus_{i\in \mathbb{Z}}H_i(X)$. Then we show that
$$\depth X-\cd_{\fa}(H(X)^\sharp)\leq \fdepth(\fa,X)\leq
\dim X-\cd_{\fa}(H(X)^\sharp).$$
This immediately provides a new characterization of Cohen-Macaulay
complexes, see Corollary 3.7 below.

\section{Formal Local Cohomology and Local Tate Cohomology}

Throughout, the symbol $\simeq$ will denote isomorphisms in the
category $\mathcal{D}(R)$, the derived category of $R$-modules. We
denote the full subcategory of $R$-modules by $\mathcal{C}_0(R)$.
The full subcategory of homologically bounded complexes is denoted
by $\mathcal{D}_{\Box}(R)$ and that of complexes homologically
bounded to the right (resp. left) is denoted by
$\mathcal{D}_{\sqsupset}(R)$ (resp. $\mathcal{D}_{\sqsubset}(R)$).
Also, if $\sharp$ is one of the symbols $\Box$ or $\sqsupset$, then
$\mathcal{D}_{\sharp}^f(R)$ stands for the full subcategory of
complexes $X\in \mathcal{D}_{\sharp}(R)$ whose all homology modules
are finitely generated.

Let $\fa$ be an ideal of $R$. The right derived functor of the
$\fa$-section functor
$$\Gamma_{\fa}(-):=\underset{n}{\varinjlim}\Hom_R(R/\fa^n,-)
:\mathcal{C}_0(R)\lo
\mathcal{C}_0(R)$$ exists in $\mathcal{D}(R)$. Let
$X\in\mathcal{D}(R)$. Then the complex ${\bf R}\Gamma_{\fa}(X)$ is
defined by ${\bf R}\Gamma_{\fa}(X) :=\Gamma_{\fa}(I)$, where $I$ is
any $K$-injective resolution of $X$. (For more details on the theory
of $K$-resolutions, we refer the reader to \cite{Sp}.) For any integer
$i$, the $i$-th
local cohomology module of $X$ with respect to $\fa$ is defined by
$H^i_{\fa}(X):=H_{-i}({\bf R}\Gamma_{\fa}(X))$. Let
$\Check{C}(\underline{\fa})$ denote the \v{C}ech complex of
$R$ on a set $\underline{\fa}=a_1,\ldots, a_n$ of generators of $\fa$. By
\cite[Theorem 1.1 iv)]{Sc1}, $${\bf R}\Gamma_{\fa}(X)\simeq
X\otimes_R^{{\bf L}} \Check{C}(\underline{\fa}). \   \  (*)$$ Let
$X\in \mathcal{D}^f_{\sqsupset}(R)$ and  $Y\in
\mathcal{D}_{\Box}(R)$. As $\Check{C}(\underline{\fa})$ is a bounded
complex of flat $R$-modules, tensor evaluation property (see
\cite[A.4.23]{C}) along with $(*)$ yield that $${\bf
R}\Gamma_{\fa}({\bf R}\Hom_R(X,Y))={\bf R}\Hom_R(X,{\bf
R}\Gamma_{\fa}(Y)). \  \  (**)$$

The left derived functor of the $\fa$-adic completion functor
$$\Lambda^{\fa}(-)=\underset{n}{\vpl}(R/\fa^n\otimes_R-):
\mathcal{C}_0(R)\lo \mathcal{C}_0(R)$$ exists in $\mathcal{D}(R)$,
and so for a complex $X\in \mathcal{D}(R)$, the complex ${\bf
L}\Lambda^{\fa}(X)$ is defined by ${\bf
L}\Lambda^{\fa}(X):=\Lambda^{\fa}(P)$, where $P$ is any
$K$-projective resolution of $X$. For any integer $i$, the $i$-th
local homology module of a complex $X\in \mathcal{D}(R)$ with
respect to $\fa$ is defined by $H_i^{\fa}(X):=H_i({\bf
L}\Lambda^{\fa}(X)).$ By \cite[(0.3), aff, page 4]{LLT} (see also
\cite[Section 4]{Sc1} for corrections) for any $X\in
\mathcal{D}(R)$, one has
$${\bf L}\Lambda^{\fa}(X)\simeq {\bf R}\Hom_R(\Check{C}(\underline{\fa}),X).
\    \   (\dag)$$ Using adjointness, $(\dag)$ and $(*)$ yields the
following isomorphisms
$${\bf L}\Lambda^{\fa}({\bf R}\Hom_R(X,Y))\simeq {\bf R}\Hom_R({\bf R}
\Gamma_{\fa}(X),Y)\simeq {\bf R}\Hom_R(X,{\bf L}\Lambda^{\fa}(Y))\
(\ddag)$$ in $\mathcal{D}(R)$ for all complexes $X$ and $Y$, see
e.g. \cite[2.6]{Fr}. In the sequel, we will use the isomorphisms
$(*), (**), (\dag)$ and $(\ddag)$ without any further comments.

\begin{definition} Let $\fa,\fb$ be two ideals of $R$ and $X\in
\mathcal{D}(R)$.
We define {\it formal local complex} of $X$ with respect to
$(\fa,\fb)$ by $\mathfrak{F}_{\fa,\fb}(X):={\bf L}\Lambda^{\fa}({\bf R}
\Gamma_{\fb}(X))\in \mathcal{D}(R)$. Also,
for each integer $i$, $i$-th {\it formal local cohomology module} of $X$
with respect to $(\fa,\fb)$ is defined by
$\mathfrak{F}_{\fa,\fb}^i(X):=H_{-i}(\mathfrak{F}_{\fa,\fb}(X))$.
When $R$ is local and $\fb$ is its maximal ideal, we abbreviate
$\mathfrak{F}_{\fa,\fb}(X)$ and $\mathfrak{F}_{\fa,\fb}^i(X)$,
respectively, by  $\mathfrak{F}_{\fa}(X)$ and
$\mathfrak{F}_{\fa}^i(X)$.
\end{definition}

A complex $X\in \mathcal{D}(R)$ is said to have finite injective dimension if
it is isomorphic, in $\mathcal{D}(R)$, to a bounded complex of injective $R$-modules.
A dualizing complex of $R$ is a complex $D\in \mathcal{D}_{\Box}^f(R)$ such that the
homothety morphism $R\longrightarrow {\bf R}\Hom_R(D,D)$ is an isomorphism in
$\mathcal{D}(R)$ and $D$ has finite injective dimension. A dualizing
complex $D$ is said to be normalized if $\sup D=\dim R$. Assume that $R$
possesses a dualizing complex $D$. Then, we know that  $\dim R$
should be finite. Now, one can check easily that $\Sigma^{\dim R-\sup X}D$ is a
normalized dualizing complex of $R$. In what follows, whenever $R$ possesses
a normalized dualizing complex $D$, we denote ${\bf R}\Hom_R(-,D)$ by $(-)^\dagger$.

\begin{lemma} (Duality) Let $\fa,\fb$ be two ideals of $R$ and
$X\in \mathcal{D}^f_{\Box}(R)$. Assume that $R$ possesses a
normalized dualizing complex $D$. Then
$\mathfrak{F}_{\fa,\fb}(X)\simeq  {\bf R}\Hom_R({\bf
R}\Gamma_{\fa}(X^\dagger),{\bf R}\Gamma_{\fb}(D)).$
\end{lemma}

\begin{prf} One can see easily that there is a natural isomorphism
$(X^\dagger)^\dagger\simeq X$. Since $X^\dagger\in \mathcal{D}^f_{\Box}(R)$,
we have
$$\begin{array}{ll}
\mathfrak{F}_{\fa,\fb}(X)
&\simeq {\bf R}\Hom_R(\Check{C}(\underline{\fa}),{\bf R}\Gamma_{\fb}(X))\\
&\simeq {\bf R}\Hom_R(\Check{C}(\underline{\fa}),{\bf
R}\Gamma_{\fb}({\bf R}
\Hom_R(X^\dagger,D)))\\
&\simeq {\bf R}\Hom_R(\Check{C}(\underline{\fa}),{\bf
R}\Hom_R(X^\dagger,
{\bf R}\Gamma_{\fb}(D)))\\
&\simeq {\bf R}\Hom_R({\bf R}\Gamma_{\fa}(X^\dagger),{\bf
R}\Gamma_{\fb}(D)).
\end{array}
$$
\end{prf}

For a local ring $(R,\fm,k)$, we denote the Matlis duality functor
$\Hom_R(-,E_R(k))$  by $(-)^\vee$.

\begin{corollary} Let $(R,\fm)$ be a local ring possessing a normalized
dualizing complex $D$, $\fa$ an ideal of $R$ and $X\in
\mathcal{D}^f_{\Box}(R)$. Then $\mathfrak{F}_{\fa}(X)\simeq {\bf R}
\Gamma_{\fa}(X^\dagger)^\vee.$
\end{corollary}

\begin{prf} By \cite[Proposition 6.1]{H}, ${\bf R}\Gamma_{\fm}(D)=E_R(k)$.
Applying Lemma 2.2 to $\fb:=\fm$ yields that
$$\mathfrak{F}_{\fa}(X)\simeq {\bf R}\Hom_R({\bf
R}\Gamma_{\fa}(X^\dagger),{\bf R}\Gamma_{\fm}(D))\simeq {\bf R}
\Gamma_{\fa}(X^\dagger)^\vee .$$
\end{prf}

Next, we present the following corollary. It shows that Definition
2.1 extends Schenzel's definition.

\begin{corollary} Let $\fa$ be an ideal of a local ring $(R,\fm)$
and $X\in \mathcal{D}(R)$. If either $X$ is a bounded complex of flat
$R$-modules whose all homology modules are finitely generated or  $X$ is a
finitely generated $R$-module, then
$\mathfrak{F}_{\fa}^i(X)\cong \underset{n}{\vpl}H^i_{\fm}(X/\fa^nX)$
for all $i\in\mathbb{Z}$.
\end{corollary}

\begin{prf} First assume that $X$ is a bounded complex of flat
$R$-modules whose all homology modules are finitely generated. Then
$X\otimes_R^{{\bf L}} \Check{C}(\underline{\fm})\simeq X\otimes_R\Check{C}
(\underline{\fm})$ is a bounded complex of flat $R$-modules and
${\bf R}\Gamma_{\fm}(X)\simeq X\otimes_R\Check{C}(\underline{\fm})$.
Hence $X\otimes_R\Check{C}(\underline{\fm})$ is a $K$-flat resolution of
${\bf R}\Gamma_{\fm}(X)$.
It is known that for a complex $Z$ and any $K$-flat resolution $F$ of $Z$,
one has ${\bf L}\Lambda^{\fa}(Z)\simeq \Lambda^{\fa}(F)$. Thus
$$\mathfrak{F}_{\fa}(X)={\bf L}\Lambda^{\fa}({\bf R}\Gamma_{\fm}(X))\simeq \Lambda^{\fa}(X\otimes_R\Check{C}(\underline{\fm})).$$ Now, we have
$$\Lambda^{\fa}(X\otimes_R\Check{C}(\underline{\fm}))=
\underset{n}{\vpl}(R/\fa^n\otimes_R(X\otimes_R\Check{C}(\underline{\fm})))
\simeq \underset{n}{\vpl}(X/\fa^nX\otimes_R\Check{C}(\underline{\fm})).$$
For each nonnegative integer $n$, all homology modules of the complex
$X/\fa^nX\otimes_R\Check{C}(\underline{\fm})(\simeq {\bf
R}\Gamma_{\fm}(X/\fa^nX))$ are Artinian. The Mittag-Leffler condition
\cite[Proposition 3.5.7]{W} implies that
${\vpl}^1_nH_{-i}(X/\fa^n\otimes_R\Check{C}(\underline{\fm}))=0$ for
all integers $i$.  From \cite[Theorem 3.5.8]{W}, we deduce the exact
sequence
$$0\lo \underset{n}{\vpl}^1H_{-i+1}(X/\fa^nX\otimes_R\Check{C}(\underline{\fm}))
\lo H_{-i}(\underset{n}{\vpl}(X/\fa^nX\otimes_R\Check{C}(\underline{\fm})))\lo
\underset{n}{\vpl}H_{-i}(X/\fa^nX\otimes_R\Check{C}(\underline{\fm}))\lo 0,$$ for all $i$.
Thus for any integer $i$, we have
$$\begin{array}{ll}\mathfrak{F}_{\fa}^i(X)&\cong H_{-i}(\underset{n}{\vpl}(X/\fa^nX\otimes_R\Check{C}(\underline{\fm})))\\
&\cong \underset{n}{\vpl}H_{-i}(X/\fa^nX\otimes_R\Check{C}(\underline{\fm}))\\
&\cong \underset{n}{\vpl}H_{-i}({\bf R}\Gamma_{\fm}(X/\fa^nX))\\
&=\underset{n}{\vpl}H^i_{\fm}(X/\fa^n
X).
\end{array}
$$

Next, assume that $X$ is a finitely generated $R$-module. Without
loss of generality, we may and do assume that $R$ is complete.
Hence, $R$ possesses a normalized dualizing complex, and so
Corollary 2.3 implies that $H_{-i}(\mathfrak{F}_{\fa}(X))\cong
H^{-i}_{\fa} (X^\dagger)^\vee.$ Now, \cite[Theorem 3.5]{Sc2}
finishes the proof in this case.
\end{prf}

Next, we bring two more results concerning the computation of formal
local cohomology modules. The first one indicates that the theory of
formal local cohomology can be considered as a unified
generalization of the two theories of local cohomology and local
homology.

\begin{proposition} Let $\fa,\fb$ be two ideals of $R$ and $X\in
\mathcal{D}(R)$. The following assertions hold.
\begin{enumerate}
\item[i)] $\mathfrak{F}_{\fa,\fa}(X)\simeq {\bf L}\Lambda^{\fa}(X)$.
\item[ii)] Assume that $\Supp_RX\subseteq \V(\fb)$. Then $\mathfrak{F}_{\fa,\fb}(X)
\simeq {\bf L}\Lambda^{\fa}(X)$.
\item[iii)] Assume that $R$ possesses a normalized dualizing complex $D$,
$X\in \mathcal{D}^f_{\Box}(R)$ and $\Supp_RX\subseteq \V(\fa)$. Then
$\mathfrak{F}_{\fa,\fb}(X)\simeq {\bf R}\Gamma_{\fb}(X)$.
\end{enumerate}
\end{proposition}

\begin{prf}  i) holds by \cite [Cor. to (0.3)*]{LLT}.\\
ii) holds by \cite[Corollary 3.2.1]{L}.\\
iii)  For any prime ideal $\fp$ of $R$, we have
$(X^\dagger)_{\fp}\simeq {\bf R}\Hom_{R_{\fp}}(X_{\fp},D_{\fp}),$
and so $\Supp_RX^\dagger\subseteq \Supp_RX.$ Thus $X^\dagger$ is
homologically bounded and $\Supp_RX^\dagger \subseteq \V(\fa)$. Now,
\cite[Corollary 3.2.1]{L} yields that ${\bf R}\Gamma_{\fa}(X^\dagger)\simeq
X^\dagger$. Set $(-)^*:={\bf R}\Hom_R(-,{\bf R}\Gamma_{\fb}(D)).$
Then Lemma 2.2 yields that
$$\mathfrak{F}_{\fa,\fb}(X)\simeq {\bf
R}\Gamma_{\fa}(X^\dagger)^*\simeq (X^\dagger)^* \simeq {\bf
R}\Gamma_{\fb}({\bf R}\Hom_R(X^\dagger,D))\simeq{\bf
R}\Gamma_{\fb}(X).$$
\end{prf}

\begin{proposition} Let $\fa,\fb$ be two ideals of $R$,
$X\in \mathcal{D}^f_{\sqsupset}(R)$ and $Y\in \mathcal{D}_{\square}(R)$. Then
\begin{enumerate}
\item[i)] $\mathfrak{F}_{\fa,\fb}({\bf R}\Hom_R(X,Y))\simeq {\bf R}
\Hom_R(X,\mathfrak{F}_{\fa,\fb}(Y))$.
\item[ii)] $\mathfrak{F}_{\fa,\fb}(X\otimes^{{\bf L}}_RY)\simeq X
\otimes^{{\bf L}}_R\mathfrak{F}_{\fa,\fb}(Y)$.
\end{enumerate}
\end{proposition}

\begin{prf} One has
$$\begin{array}{ll}
\mathfrak{F}_{\fa,\fb}({\bf R}\Hom_R(X,Y))&\simeq {\bf
L}\Lambda^{\fa}({\bf R}
\Gamma_{\fb}({\bf R}\Hom_R(X,Y)))\\
&\simeq {\bf R}\Hom_R(\Check{C}(\underline{\fa}),{\bf R}
\Gamma_{\fb}({\bf R}\Hom_R(X,Y)))\\
&\simeq {\bf R}\Hom_R(\Check{C}(\underline{\fa}),{\bf R}
\Hom_R(X,{\bf R}\Gamma_{\fb}(Y)))\\
&\simeq {\bf R}\Hom_R(X,{\bf R}\Hom_R(\Check{C}(\underline{\fa}),
{\bf R}\Gamma_{\fb}(Y)))\\
&\simeq{\bf R}\Hom_R(X,\mathfrak{F}_{\fa,\fb}(Y)).
\end{array}
$$

ii) One has
$$\begin{array}{ll} \mathfrak{F}_{\fa,\fb}(X\otimes^{{\bf L}}_RY)
&\simeq {\bf L}\Lambda^{\fa}({\bf R}
\Gamma_{\fb}(X\otimes^{{\bf L}}_RY))\\
&\simeq {\bf R}\Hom_R(\Check{C}(\underline{\fa}),{\bf R}
\Gamma_{\fb}(X\otimes^{{\bf L}}_RY))\\
&\simeq {\bf R}\Hom_R(\Check{C}(\underline{\fa}),{\bf R}
\Gamma_{\fb}(Y)\otimes^{{\bf L}}_RX)\\
&\overset{\sharp}\simeq {\bf
R}\Hom_R(\Check{C}(\underline{\fa}),{\bf R}
\Gamma_{\fb}(Y))\otimes^{{\bf L}}_RX\\
&\simeq X\otimes^{{\bf L}}_R\mathfrak{F}_{\fa,\fb}(Y).
\end{array}
$$
Regarding the isomorphism $\sharp$, we have to give some
explanations. By \cite[5.8]{CFH}, the projective dimension of
$\Check{C}(\underline{\fa})$ is finite. On the other hand, as $Y\in
\mathcal{D}_{\square}(R)$ and $\Check{C}(\underline{\fa})$ is a
bounded complex of flat $R$-modules, we deduce that ${\bf
R}\Gamma_{\fa}(Y)\simeq \Check{C}(\underline{\fa})\otimes^{{\bf
L}}_RY$ is homologically bounded. Thus \cite[Proposition 2.2
vi)]{CH} implies the isomorphism $\sharp$.
\end{prf}

The theory of local Tate cohomology was introduced by Greenlees
\cite{G}. Let $\{\phi^t:X^t\lo X^{t+1}\}_{t\in \mathbb{N}}$ be a
family of morphisms of complexes. It induces a morphism of complexes
$\varphi=(\varphi_i):\bigoplus_{t\in \mathbb{N}}X^t\lo
\bigoplus_{t\in \mathbb{N}}X^t$ given by
$\varphi_i((x_i^t))=(x_i^t)-(\phi^t_i(x_i^t))$ in spot $i$. The
telescope of $\{\phi^t:X^t\lo X^{t+1}\}_{t\in \mathbb{N}}$ is
defined by $\Tel(X^t):=\Cone(\varphi)$. Let $a\in R$. For each
natural integer $t$, let $K(a^t)$ denote the Koszul complex of $R$
with respect to $a^t$. Clearly, multiplication by $a$, induces a
family of morphisms of complexes $K(a^t)\lo K(a^{t+1})$. The
projective stabilized Koszul complex with respect to $a$ is defined
by $K(a^{\infty}):=\Tel (K(a^t))$. Let $\underline{\fa}=a_1,\ldots,
a_n$ be a sequence of elements of $R$. The projective stabilized
Koszul complex with respect to $\underline{\fa}$ is defined by
$K(\underline{\fa}^{\infty}):=\Tel (K(a_1^t))\otimes_R\ldots
\otimes_R \Tel (K(a_n^t))$. The stabilized  \v{C}ech complex
with respect to $\underline{\fa}$ is defined by
$\Check{C}(\underline{\fa}^{\infty}):=
\Cone(K(\underline{\fa}^{\infty})\lo R)$.

\begin{definition} Let $\fa$ be an ideal of $R$ and $X\in \mathcal{D}(R)$.
Let $\underline{\fa}=a_1,\ldots, a_n$ be a generating set of $\fa$.
The {\it local Tate complex} of $X$ with respect to $\fa$ is defined
by $T(X):=
\Hom_R(K(\underline{\fa}^{\infty}),X)\otimes_R\Check{C}(\underline{\fa}^{\infty})$.
Also, for each integer $i$, $i$-th {\it local Tate cohomology
module} of $X$ with respect to $\fa$ is defined by
$\widehat{H}^i_{\fa}(X):=H_{-i}(T(X))$.
\end{definition}

\begin{lemma} (The algebraic Warwick duality) Let $\underline{\fa}=a_1,\ldots, a_n$ be a
sequence of elements of $R$ and $X\in \mathcal{D}(R)$. Let
$TT(X):=\Hom_R(\Check{C}(\underline{\fa}^{\infty}), X\otimes_R\Sigma
K(\underline{\fa}^{\infty}))$. Then there is a natural isomorphism
$T(X)\simeq TT(X)$ in $\mathcal{D}(R)$.
\end{lemma}

\begin{prf} For modules, this is proved in \cite[Theorem 4.1]{G}.
In view of \cite[Corollary 4.6]{G}, we can check easily that the
conclusion is true also for complexes.
\end{prf}

The next result is contained in \cite[Proposition 5.1.3]{LLT}. Here
we prove it by applying a more simple argument.

\begin{proposition} Let $\fa$ be an ideal of $R$ and $X\in
\mathcal{D}_{\Box}(R)$. We have the following long exact sequence
$$\cdots \lo H_{\fa}^i(X)\lo H^{\fa}_{-i}(X)\lo
\widehat{H}^i_{\fa}(X)\lo H_{\fa}^{i+1}(X)\lo H^{\fa}_{-i-1}(X)\lo
\cdots .$$
\end{proposition}

\begin{prf} Let $\underline{\fa}=a_1,\ldots, a_n$ be a generating set of $\fa$.
In view of Lemma
2.8, we have the following natural isomorphisms in $\mathcal{D}(R)$:
$$\begin{array}{ll} T(X)&\simeq TT(X)\\
&\simeq \Sigma\Hom_R(\Check{C}(\underline{\fa}^{\infty}),
X\otimes_RK(\underline{\fa}^{\infty}))\\
&=\Sigma\Hom_R(\Cone(K(\underline{\fa}^{\infty})\lo R),
X\otimes_RK(\underline{\fa}^{\infty}))\\
&\simeq \Sigma
\Sigma^{-1}\Cone(\Hom_R(R,X\otimes_RK(\underline{\fa}^{\infty})) \lo
\Hom_R(K(\underline{\fa}^{\infty}),X\otimes_RK
(\underline{\fa}^{\infty})))\\
&\simeq \Sigma
\Sigma^{-1}\Cone(X\otimes_RK(\underline{\fa}^{\infty}) \lo
\Hom_R(K(\underline{\fa}^{\infty}),X\otimes_RK
(\underline{\fa}^{\infty})))\\
&=\Cone(X\otimes_RK(\underline{\fa}^{\infty})\lo
\Hom_R(K(\underline{\fa}^{\infty}),X\otimes_RK
(\underline{\fa}^{\infty}))).
\end{array}
$$
Note that the third quism above is obtained by \cite[3.48]{Fo2}. In
view of \cite[Remark 1.1]{G}, $K(\underline{\fa}^{\infty})$ is a
projective resolution of $\Check{C}(\underline{\fa})$. So, we deduce
the following exact sequence in $\mathcal{D}(R)$
$$0\lo {\bf R}\Hom_R(\Check{C}(\underline{\fa}),
X\otimes^{{\bf L}}_R\Check{C}(\underline{\fa}))\lo T(X)\lo \Sigma
(X\otimes^{{\bf L}}_R\Check{C}(\underline{\fa}))\lo 0. $$ Now, by
Proposition 2.5 i), the induced long exact sequence of homologies of
this exact sequence is our desired long exact sequence.
\end{prf}

Let $\fa$ be an ideal of $R$ and $X\in \mathcal{D}(R)$. We denote
$$\sup\{i\in \mathbb{Z}|H^i_{\fa}(X)\neq 0\}(=-\inf {\bf
R}\Gamma_{\fa}(X))$$ by $\cd_{\fa}(X)$. Also, recall that
$\depth(\fa,X)$ is defined by $\depth(\fa,X):=-\sup {\bf
R}\Hom_R(R/\fa,X)$. The following corollary is immediate.

\begin{corollary} Let $\fa$ be an ideal of $R$ and
$X\in \mathcal{D}_{\Box}(R)$. If either $i<\depth(\fa,X)-1$ or
$i>\cd_{\fa}(X)$, then $\widehat{H}^i_{\fa}(X)\cong
H^{\fa}_{-i}(X)$.
\end{corollary}

\begin{corollary} Let $\fa$ be an ideal of $R$ and $M$ an $R$-module.
Then $\widehat{H}_{\fa}^i(M)\cong H_{\fa}^{i+1}(M)$ for all $i>0$.
Moreover, if $M$ is $\fa$-adic complete, then
$\widehat{H}_{\fa}^i(M)\cong H_{\fa}^{i+1}(M)$ for all $i\neq 0,-1$.
\end{corollary}

\begin{prf} Since $H_{-i}^{\fa}(M)=0$ for all $i>0$, from Proposition 2.9, it
turns out that $\widehat{H}_{\fa}^i(M)\cong H_{\fa}^{i+1}(M)$ for
all $i>0$. Now, assume that $M$ is $\fa$-adic complete. Then
$H^{\fa}_i(M)=0$ for all $i\neq 0$. Hence from Proposition 2.9, we
deduce that $\widehat{H}_{\fa}^i(M)\cong H_{\fa}^{i+1}(M)$ for all
$i\neq 0,-1$.
\end{prf}

\section{Vanishing Results}

We start this section with the following definition.

\begin{definition} Let $\fa,\fb$ be two ideals of $R$ and
$X\in \mathcal{D}(R)$. We define {\it formal depth} of $X$ with
respect to $(\fa,\fb)$ by $\fdepth(\fa,\fb,X):=\inf\{i\in
\mathbb{Z}:\mathfrak{F}_{\fa,\fb}^i(X)\neq 0\}.$ When $R$ is local
with maximal ideal $\fm$, we abbreviate $\fdepth(\fa,\fm,X)$ by
$\fdepth(\fa,X)$.
\end{definition}

Let $\fa$ be an ideal of $R$ and $X\in \mathcal{D}(R)$. Recall that
$\dim_RX$ is defined by $\dim_RX:=\sup\{\dim R/\fp-\inf
X_{\fp}|\fp\in \Spec R\}$. It is known that for any complex $X\in
\mathcal{D}_{\Box}^f(R)$, we have $-\sup {\bf
R}\Gamma_{\fa}(X)=\depth(\fa,X)$ and $\cd_{\fa}(X)\leq \dim_RX$ with
equality if $R$ is local and $\fa$ is its maximal ideal, see
\cite[Theorem 7.8 and Proposition 7.10]{Fo1}. From Corollary 2.3, we
can record the following immediate corollary.

\begin{corollary} Let $(R,\fm)$ be a local ring possessing a normalized
dualizing complex $D$, $\fa$ an ideal of $R$ and $X\in
\mathcal{D}^f_{\Box}(R)$. Then
$$\fdepth(\fa,X)=-\sup\{i\in\mathbb{Z}:H^{\fa}_i({\bf R}\Gamma_{\fm}(X))
\neq 0\}=-\cd_{\fa}(X^\dagger)\geq
-\dim_RX^\dagger.$$
\end{corollary}

Next, is our first main result.

\begin{theorem} Let $\fa,\fb$ be two ideals of $R$ and $X\in \mathcal{D}^f_{\Box}(R)$.
Let $K(\underline{\fa})$ denote the Koszul complex of $R$ with respect to a generating set
$\underline{\fa}=a_1,\ldots, a_n$ of $\fa$. The following assertions hold.
\begin{enumerate}
\item[i)] $\fdepth(\fa,\fa,X)=-\sup {\bf L}\Lambda^{\fa}(X)$.
\item[ii)] $\fdepth(\fa,\fb,X)\geq \depth (\fb,X)-\cd_{\fa}(R).$
\item[iii)]  $\sup\{i\in\mathbb{Z}:\mathfrak{F}_{\fa,\fb}^i(X)\neq 0\}=
\cd_{\fb}(K(\underline{\fa})\otimes_R^{\bf L}X)\leq \dim_RX/\fa X$. In particular,
if $R$ is local, then  $\sup\{i\in\mathbb{Z}:\mathfrak{F}_{\fa}^i(X)\neq 0\}=\dim_RX/\fa X$.
\end{enumerate}
\end{theorem}

\begin{prf} i) is clear by Proposition 2.5 i).

ii) For any two complexes $V\in \mathcal{D}_{\sqsupset}(R)$ and
$W\in \mathcal{D}_{\sqsubset}(R)$, \cite[Proposition A.4.6]{C} yields that
$$\sup {\bf R}\Hom_R(V,W)\leq \sup W-\inf V.$$ Hence,
one has
$$\begin{array}{ll}
\fdepth(\fa,\fb,X)&=\inf\{i\in\mathbb{Z}:H_{-i}({\bf L}\Lambda^{\fa}
({\bf R}\Gamma_{\fb}(X)))\neq 0\}\\
&=-\sup\{i\in\mathbb{Z}:H_{i}({\bf L}\Lambda^{\fa}({\bf R}
\Gamma_{\fb}(X)))\neq 0\}\\
&=-\sup {\bf R}\Hom_R(\Check{C}(\underline{\fa}),
{\bf R}\Gamma_{\fb}(X))\\
&\geq-\sup {\bf R}\Gamma_{\fb}(X)+\inf \Check{C}(\underline{\fa})\\
&=\depth(\fb,X)-\cd_{\fa}(R).
\end{array}
$$

iii) For any complex $Y\in \mathcal{D}_{\sqsupset}(R)$, \cite[Theorem 2.11]{Fr}
asserts that $\inf {\bf L}\Lambda^{\fa}(Y)=\inf (K(\underline{\fa})\otimes_R^{\bf L}Y)$.
Also, by
Grothendieck's vanishing Theorem (see \cite[Theorem 7.8]{Fo2}), for any complex
$Z\in \mathcal{D}^f_{\Box}(R)$, we know that $\cd_{\fb}(Z)\leq \dim_RZ$ with
equality if $R$ is local and $\fb$ is its maximal ideal. One has
$$\begin{array}{ll}
\sup\{i\in\mathbb{Z}:\mathfrak{F}_{\fa,\fb}^i(X)\neq 0\}&=\sup\{i\in\mathbb{Z}:
H_{-i}({\bf L}\Lambda^{\fa}
({\bf R}\Gamma_{\fb}(X)))\neq 0\}\\&=-\inf {\bf L}\Lambda^{\fa}({\bf R}
\Gamma_{\fb}(X))\\
&=-\inf (K(\underline{\fa})\otimes_R^{\bf L}{\bf R}\Gamma_{\fb}(X))\\
&=-\inf {\bf R}\Gamma_{\fb}(K(\underline{\fa})\otimes_R^{\bf L}X)\\
&=\cd_{\fb}(K(\underline{\fa})\otimes_R^{\bf L}X)\\
&\leq \dim_R(K(\underline{\fa})\otimes_R^{\bf L}X)\\
&=\dim_RX/\fa X.
\end{array}
$$
It is easy to see that $\Supp_R(K(\underline{\fa})\otimes_R^{\bf L}X)=\Supp_RX/\fa X$
and for any prime ideal $\fp$ of $R$, we have $\inf(K(\underline{\fa})
\otimes_R^{\bf L}X)_{\fp}=\inf X_{\fp}$.
This yields the last equality above.
\end{prf}

Let $\fa$ be an ideals of $R$ and $X\in \mathcal{D}^f_{\Box}(R)$.
Finding a good upper bound for $\sup {\bf L}\Lambda^{\fa}(X)$ is of some
interest (see e.g. \cite[page 179]{Sc1}).  Theorem 3.3
immediately implies the following corollary.

\begin{corollary} Let $\fa$ be an ideal of $R$ and $X\in
\mathcal{D}^f_{\Box}(R)$. Then
$\sup {\bf L}\Lambda^{\fa}(X)\leq \cd_{\fa}(R)-\depth(\fa,X).$
\end{corollary}

For proving our second main result, we need to prove the following lemma.

\begin{lemma} Let $\fa$ be an ideal of $R$ and $M$ a finitely generated $R$-module.
Let $\{N_\lambda\}_{\lambda\in \Lambda}$ be a family of finitely
generated $R$-modules such that $\cup_{\lambda\in
\Lambda}\Supp_RN_\lambda=\Supp_RM$. Then
$$\cd_{\fa}(M)=\sup\{\cd_{\fa}(N_\lambda)|\lambda\in \Lambda\}.$$
\end{lemma}

\begin{prf} By \cite[Theorem 2.2]{DNT}, for any finitely generated $R$-module
$N$, we conclude that $\cd_{\fa}(N)=\cd_{\fa}(\oplus_{\fp\in \Supp_RN}R/\fp)$. Hence,
$$\cd_{\fa}(N)=\sup\{\cd_{\fa}(R/\fp)|\fp\in \Supp_RN\},$$
which easily yields the claim.
\end{prf}

The following result extends \cite[Theorem 4.6]{AD} to complexes. It
should be noted that its proof is completely different from our
proof for \cite[Theorem 4.6]{AD}.

\begin{theorem} Let  $\fa$ be an ideal of a local ring $(R,\fm)$ and
$X\in \mathcal{D}^f_{\Box}(R)$ a non-homologically trivial complex. Set
$H(X)^\sharp:=\oplus_{i\in
\mathbb{Z}}H_i(X)$. Then
$$\depth X-\cd_{\fa}(H(X)^\sharp)\leq \fdepth(\fa,X)\leq \dim X-
\cd_{\fa}(H(X)^\sharp).$$
\end{theorem}

\begin{prf} \cite[Corollary 3.4.4]{L} yields an isomorphism
$${\bf R}\Gamma_{\fm}(X) \otimes_R\hat{R}\simeq {\bf R}\Gamma_{\fm\hat{R}}
(X\otimes_R\hat{R})$$ in $\mathcal{D}(R)$. Let $P$ be a $K$-projective
resolution of ${\bf R}\Gamma_{\fm}(X)$. Then $P\otimes_R\hat{R}$ is a
$K$-projective resolution of the $\hat{R}$-complex ${\bf R}\Gamma_{\fm\hat{R}}
(X\otimes_R\hat{R})$. Now, consider the following natural isomorphisms
of $R$-complexes
$$\begin{array}{ll}\Lambda^{\fa \hat{R}}(P\otimes_R\hat{R})&=
\underset{n}\vpl (\hat{R}/(\fa \hat{R})^n\otimes_{\hat{R}}
(P\otimes_R\hat{R}))\\
&\simeq \underset{n}\vpl (R/\fa^n\otimes_RP)\\
&=\Lambda^{\fa}(P),
\end{array}$$
which implies that the two complexes $\mathfrak{F}_{\fa}(X)$ and
$\mathfrak{F}_{\fa \hat{R}}(X\otimes_R\hat{R})$ are isomorphic in
$\mathcal{D}(R)$. Thus $\fdepth(\fa,X)=\fdepth(\fa \hat{R},X\otimes_R\hat{R})$.
Next, it is  straightforward to check that $\depth X=\depth (X\otimes_R\hat{R})$
and $\cd_{\fa}(H(X)^\sharp)=\cd_{\fa \hat{R}}(H(X\otimes_R\hat{R})^\sharp)$.
On the other hand, by \cite[Proposition 3.5]{Fo3}, one has
$$\begin{array}{ll} \dim X&=\sup\{\dim H_i(X)-i|i\in \mathbb{Z}\}\\
&=\sup\{\dim H_i(X\otimes_R\hat{R})-i|i\in \mathbb{Z}\}\\
&=\dim (X\otimes_R\hat{R}).
\end{array}$$
Therefore, without loss of generality,
we may and do assume that $R$ is complete. So, $R$ possesses a normalized
dualizing complex $D$. By Corollary 3.2,
$\fdepth(\fa,X)=-\cd_{\fa}(X^\dagger)$.  By \cite[16.20]{Fo2}, it
turns out that $\inf X^\dagger=\depth X$ and $\sup X^\dagger=\dim
X$. The natural quism $(X^\dagger)^\dagger\simeq X$ yields that
$\Supp_RX^\dagger=\Supp_RX$, and so
$$\bigcup_{i\in \mathbb{Z}}\Supp_RH_i(X^\dagger)=\Supp_RH(X)^\sharp.$$
Hence Lemma 3.5 implies that
$$\cd_{\fa}(H(X)^\sharp)=\sup\{\cd_{\fa}(H_i(X^\dagger))|i\in
\mathbb{Z}\}.$$ In particular, there exists $\depth X\leq j\leq \dim
X$ such that $\cd_{\fa}(H(X)^\sharp)=\cd_{\fa}(H_j(X^\dagger))$. Since
$X^\dagger\in \mathcal{D}^f_{\Box}(R)$, by \cite[Theorem 3.2]{DY}, one has
$$\cd_{\fa}(X^\dagger)=\sup\{\cd_{\fa}(H_i(X^\dagger))-i|\depth X\leq
i\leq \dim X \}.$$
Thus $$\cd_{\fa}(H(X)^\sharp)-\dim X\leq
\cd_{\fa}(H_j(X^\dagger))-j\leq \cd_{\fa}(X^\dagger)\leq
\cd_{\fa}(H(X)^\sharp)-\depth X,$$
which is what we wish to prove it.
\end{prf}

The next result provides a characterization for Cohen-Macaulay complexes.

\begin{corollary} Let $(R,\fm)$ be a local ring and $X\in \mathcal{D}^f_{\Box}(R)$
a non-homologically trivial complex. Set $H(X)^\sharp:=\oplus_{i\in
\mathbb{Z}}H_i(X)$. The following are equivalent:
\begin{enumerate}
\item[i)] $X$ is Cohen-Macaulay.
\item[ii)] $\fdepth(\fa,X)=\dim X-\cd_{\fa}(H(X)^\sharp)$ for
all ideals $\fa$ of $R$.
\item[iii)] $\fdepth(0,X)=\dim X-\cd_{0}(H(X)^\sharp)$.

\end{enumerate}
\end{corollary}

\begin{prf} By Theorem 3.6, the implication $i)\Rightarrow ii)$ is obvious.
Also clearly, ii) implies iii).

We have $$\mathfrak{F}_{0}^i(X)=
H_{-i}({\bf L}\Lambda^{0} ({\bf R}\Gamma_{\fm}(X)))\cong H_{\fm}^i(X).$$ Hence
$\fdepth(0,X)=\depth X$. Also, we can see easily that
$\cd_{0}(H(X)^\sharp)=0$. Thus iii) yields i).
\end{prf}

\begin{remark}  Let $\fa$ be an ideal of a local ring $(R,\fm)$.
\begin{enumerate}
\item[i)] Suppose that $\fa$ is generated by a regular sequence
$\underline{x}:=
x_{1},\ldots, x_{r}$   and $K(\underline{x})$ denotes the Koszul
complex of $R$ with respect to $\underline{x}$. Then
$\fdepth(\fa,K(\underline{x}))=\depth R-\Ht \fa$. To this end, note
that $R/\fa\simeq K(\underline{x})$. Hence
$$\mathfrak{F}_{\fa}^i(K(\underline{x}))\cong\mathfrak{F}_{\fa}^i(R/\fa)
\cong
H^i_{\fm}(R/\fa)$$ for all $i\geq 0$, and so
$$\fdepth(\fa,K(\underline{x}))=\depth R/ \fa=\depth R-\Ht \fa.$$
\item[ii)] Let $X\in \mathcal{D}^f_{\Box}(R)$ be a $d$-dimensional
Cohen-Macaulay complex. Then $$\mathfrak{F}^i_{\fa}(X)\cong
H_{\fa}^{d-i}(H_{\fm}^{d}(X)^\vee)^\vee \cong
H^{\fa}_{d-i}(H_{\fm}^d(X)) .$$  To see this, first of all note that
without loss of generality, we may and do assume that $R$ is
complete. So, $R$ possesses a normalized dualizing complex. Since $X$
is Cohen-Macaulay, one has $H_{\fm}^{i}(X)=0$ for all $i\neq d$.
Hence ${\bf R}\Gamma_{\fm}(X)\simeq\Sigma^{-d}H_{\fm}^{d}(X)$.
By local duality Theorem (see \cite[Chapter V, Theorem 6.2]{H}),
one has $(X^\dagger)^\vee \simeq {\bf R}\Gamma_{\fm}(X)$. Since $R$ is complete
and $X^\dagger\in \mathcal{D}^f_{\Box}(R)$, this implies
that $X^\dagger \simeq {\bf R}\Gamma_{\fm}(X)^\vee$.
Thus by Corollary 2.3, we have
$$\mathfrak{F}_{\fa}(X)\simeq {\bf R}\Gamma_{\fa}({\bf R}
\Gamma_{\fm}(X)^\vee)^\vee \simeq {\bf
R}\Gamma_{\fa}(\Sigma^dH_{\fm}^{d}(X)^\vee)^\vee.$$ Hence
$$\mathfrak{F}^i_{\fa}(X)\simeq H_{-i}({\bf
R}\Gamma_{\fa}(\Sigma^{d}H_{\fm}^{d}(X)^\vee)^\vee)\simeq
H_{i-d}({\bf R}\Gamma_{\fa}(H_{\fm}^{d}(X)^\vee))^\vee=H_{\fa}^{d-i}
(H_{\fm}^{d}(X)^\vee)^\vee.$$ Also, we have
$$\mathfrak{F}_{\fa}(X)={\bf L}\Lambda^{\fa}({\bf
R}\Gamma_{\fm}(X))={\bf
L}\Lambda^{\fa}(\Sigma^{-d}H_{\fm}^{d}(X)),$$ and so
$$\mathfrak{F}^i_{\fa}(X)\simeq H_{-i}({\bf L}\Lambda^{\fa}
(\Sigma^{-d}H_{\fm}^{d}(X)))\simeq H_{d-i}({\bf L}\Lambda^{\fa}
(H_{\fm}^{d}(X)))=H^{\fa}_{d-i}(H_{\fm}^{d}(X)).$$
\item[iii)] Suppose that $R$ is Cohen-Macaulay and complete and
$\omega_R$ is a
canonical module of $R$. Assume that $\fa$ is cohomologically
complete intersection (i.e. $\cd_{\fa}(R)=\Ht \fa$) and let
$t:=\dim R/\fa$. Then
$\id_R(\mathfrak{F}^t_{\fa}(\omega_R))<\infty$. To this end, first
note that $\Supp_R\omega_R=\Spec R$,  and so  by \cite[Theorem
2.2]{DNT}, we have
$$\dim_R(\omega_R/\fa \omega_R)=\dim_R R/\fa=\dim R-\cd_{\fa}(R)=\dim
R-\cd_{\fa}(\omega_R)=\fdepth(\fa,\omega_R).$$ Hence
$\mathfrak{F}^i_{\fa}(\omega_R)=0$ for all $i\neq t$, and so
$$\mathfrak{F}^t_{\fa}(\omega_R)\simeq
\Sigma^t \mathfrak{F}_{\fa}(\omega_R)\simeq \Sigma^t {\bf R}
\Hom_R(\Check{C}(\underline{\fa}),\omega_R\otimes_R^{\bf
L}\Check{C}(\underline{\fm})).$$ Now, since
$\id_R(\omega_R)<\infty$, the conclusion follows by
\cite[Proposition 2.4]{CH}.
\item[iv)] The notion of Frobenius depth was defined by Hartshorne and
Speiser in \cite[page 60]{HS}. Suppose $R$ is regular of prime
characteristic. By \cite[Theorem 4.3]{Ly}, it follows that Frobenius
depth of $R/\fa$ is equal to $\dim R-\cd_{\fa}(R)$. Having
\cite[Lemma 4.8 d)]{Sc2} and Corollary 2.4 in mind, we conclude that
Frobenius depth of $R/\fa$ is equal to $\fdepth(\fa,R)$.
\item[v)] Suppose $R$ is a one-dimensional domain and set
$M:=R/\fm\oplus R$. Clearly $M$ is not Cohen-Macaulay, while
$\fdepth(\fb,M)=\dim M-\cd_{\fb}(H(M)^\sharp)$ for all non-zero
ideals $\fb$ of $R$.
\end{enumerate}
\end{remark}

\begin{acknowledgement}  We thank the anonymous referee for his/her detailed
review.
\end{acknowledgement}


\end{document}